\RequirePackage{fix-cm}
\documentclass{svjour3}                     
\smartqed  
\usepackage{graphicx}
\usepackage[top=1.1in, bottom=1.1in, left=1.2in, right=1.2in]{geometry}
\usepackage{lineno,hyperref}
\usepackage{bbm}
\usepackage{amssymb}
\usepackage[misc]{ifsym}
\usepackage{multirow,booktabs}
\usepackage{caption}
\usepackage{mathrsfs}
\usepackage{amsfonts}
\usepackage{colortbl,dcolumn}
\usepackage{tabularx}
\usepackage{amsmath}
\usepackage{verbatim}
\usepackage{booktabs}
\newcommand{\R}{\mathbb{R}}
\newcommand{\N}{\mathbb{N}}
\newcommand{\E}{\mathbb{E}}
\renewcommand{\P}{\mathbb{P}}
\newcommand{\dd}{\text{d}}

\usepackage[numbers,sort&compress]{natbib}
\begin{document}

\title{On the backward Euler method
	for a generalized Ait-Sahalia-type rate model with Poisson jumps}


\author{ Yuying Zhao   \and
       Xiaojie Wang  \and Mengchao Wang \Letter    }
\institute{ Yuying Zhao  \at
              School of Mathematics and Statistics, Central South University, Changsha, China \\
              \email{y.y.zhao68@csu.edu.cn}           
           \and
           Xiaojie Wang  \at
          School of Mathematics and Statistics, Central South University, Changsha, 
          China \\
          \email{x.j.wang7@csu.edu.cn}  
           \and 
          Mengchao Wang \Letter  \at
          School of Mathematics and Statistics, Central South University, Changsha, 
          China \\
          \email{wangmengchao@csu.edu.cn} 
          \and
           This work was supported by Natural Science Foundation of China (11671405, 11971488,
          	91630312), Innovation Program of Central
          	South University(No.2019zzts397) and 
	Natural Science Foundation of Hunan Province for Distinguished Young Scholars  (2020JJ2040). 
       }
\date{}
\maketitle

\begin{abstract}
	This article aims to reveal the mean-square convergence
	rate of the backward Euler method (BEM)
	for a generalized Ait-Sahalia interest rate model
	with Poisson jumps. The main difficulty in the analysis is caused by
	the non-globally Lipschitz drift and diffusion coefficients of the model.
	We show that the BEM preserves the positivity of
	the original problem. Furthermore, we successfully recover
	the mean-square convergence rate of order one-half
	for the BEM. The
	theoretical findings are accompanied by several numerical examples.
	\keywords{Ait-Sahalia model, Poisson jumps, backward Euler method, mean-square convergence rate \\
		AMS subject classification: {\rm\small 60H35, 60H15, 65C30.}}
\end{abstract}
\section{Introduction} \label{sect:introduction}
As is well known, stochastic differential equations (SDEs) are widely used in various scientific areas to 
model real-life phenomena affected by random noises.
However, in order to model the event-driven phenomena, it is necessary and significant to 
introduce SDEs with Poisson jumps \cite{Cont2004Financial,platen10numerical}. For instance, the stock price
movements might suffer from sudden and significant impacts
caused by unpredictable important events such as market crashes,
announcements made by central banks, changes in credit
rating, etc.
Over the last decade,  SDEs with jumps have become increasingly popular for modelling market fluctuations,
both for risk management and option pricing purposes
(see, e.g., \cite{Cont2004Financial}).
Since the analytic solutions of nonlinear SDEs with jumps are rarely available, numerical approximations 
become a powerful tool to understand the behavior of the underlying problems.
Motivated by this,  a great deal of research papers are devoted to the interesting topic
(see,
e.g., \cite{Chen2020Convergence,deng2019truncated,Li2019X,platen10numerical,Yang2017transformed,Dareiotis2016tamed,Chalmers2008tamed,Wang2010Compensated,Kaluza2018Optimal,Higham2007Strong,Buckwar2011Runge-kutta,Przybylowicz2016Optimal,Higham2005Numerical,Higham2006Convergence,Chen2019approximations,Kumar2017explicit,Kumar2017tamed}). 

The present article is concerned with the mean-square convergence analysis of a time-stepping scheme 
for a generalized Ait-Sahalia interest rate model with Poisson jumps,
which takes the form as follows,
\begin{equation}\label{eq:Ait-Sahalia-model}
\dd X_t= (a_{-1} X_t^{-1} - a_{0}+ a_{1} X_t - a_{2} X_t^{\gamma})
\, \dd t
+b X_t^{\theta} \dd W_t + \varphi(X_{t^{-}}) \, \dd N_t, \quad
X_0=x_0, 
\quad
t > 0.
\end{equation}
Here, $a_{-1}$,\,$a_{0}$,\,$a_{1}$,\,$a_{2}$,\,$ b > 0$, $\theta$,\,$\gamma>1$, 
$ \varphi : \R \rightarrow \R $ and $X_{t^{-}}:=\lim _{s  \rightarrow t^{-}} X_s $.
Moreover, we let $ \{W_t\}_{t \in [0,\infty)}$ and  $\{N_t\}_{t \in [0,\infty)}$ be a one-dimensional Brownian motion
and a Poisson process
with the jump intensity $ \lambda > 0 $, respectively, on a filtered probability space
$ (\Omega, \mathscr{ F }, \P, \{\mathscr{ F }_t\}_{t \geq 0}) $
with respect to the normal filtration $  \{\mathscr{ F }_t\}_{t \geq 0} $.
The Brownian motion and the Poisson process are assumed to be independent with each other.
By $ \tilde{N}_t := N_t-\lambda t, t \in [0,\infty) $ we denote the compensated Poisson process,
which is also a martingale with respect to the normal filtration $  \{\mathscr{ F }_t\}_{t \geq 0} $.
Recall that the generalized Ait-Sahalia interest rate model
without jumps, i.e., $ \varphi \equiv 0 $, has been already numerically studied by \cite{Szpruch2011Numerical},
where strong convergence of the backward Euler method was proved,
but without revealing any convergence rate. Taking
the jump function $ \varphi $ to be a
linear one, Deng et.al.,~\cite{Deng2018Generalized} examined
the analytical properties of the model, including the positivity,
boundedness and pathwise asymptotic
estimates. Also, they applied the Euler-Maruyama (EM) method to the particular model
and proved that the explicit scheme converges in probability to the true solution of the model.
To the best of our knowledge, no strong convergence rate has been reported in the literature for numerical approximations of
the generalized Ait-Sahalia interest rate model with Poisson jumps as \eqref{eq:Ait-Sahalia-model}.\\
\indent
In this paper, we focus on the the backward Euler method for \eqref{eq:Ait-Sahalia-model}, which has been widely studied for SDEs without jumps (see \cite{Andersson2017Mean-square,Higham2002convergence,wang19mean-square,Huang2012Exponential,Zong2014Choice,Mao2013Euler,Mao2013convergence,Neuenkirch2014First}).
When used to solve the aforementioned Ait-Sahalia model with jumps, the chosen scheme is shown to be positivity preserving.
Furthermore we provide an easy approach of error analysis to successfully obtain its mean-square convergence rate of order one-half for
full parameters in the case $ \gamma + 1 > 2 \theta $
and for parameters obeying $ \frac{a_2}{b^2} > 2\gamma - \frac{3}{2} $
in the general critical case $ \gamma + 1 = 2 \theta $.
As a byproduct, this work reveals the expected convergence rate for the model
without jumps, which is missing in \cite{Szpruch2011Numerical}.

The remainder of this paper is organized as follows.
The next section concerns the properties of the solution to the model, including the existence and uniqueness, the positivity and the moment boundedness of the solution. The mean-square convergence rate of the BEM is identified in section 3 for the generalized
Ait-Sahalia interest rate model with Poisson jumps. Finally numerical experiments are performed to illustrate the
theoretical results.

\section{The jump-extended Ait-Sahalia model} \label{sect:Ait-Sahalia}
Throughout this paper we will use the following notation.
Let $ | \cdot | $ and $ \langle \cdot ,\cdot \rangle $
be the Euclidean norm and the inner product in $ \mathbb{R} $, respectively.
Given a filtered probability space
$ ( \Omega, \mathcal{ F }, \{\mathcal{ F }_t\}_{t \in [0,T]} ,\P ) $
satisfying the usual hypotheses,
that is to say, it is right continuous and increasing while
$ \mathcal{F}_{0} $ contains all $ \P $-null sets. 
Let $ \{W_t\}_{t \in [0,\infty)} $ and  $ \{N_t\}_{t \in [0,\infty)} $ be a one-dimensional Brownian motion
and a scalar Poisson process with the jump intensity $ \lambda > 0 $, respectively,  
with respect to the normal filtration $  \{\mathscr{ F }_t\}_{t \geq 0} $.
Let $ \E $ denote the expectation and let $ L^p(\Omega;\R) $ denote the space consisting of $ \R $-value
$ p $-times integrable random variables with the norm defined by
$ \| \xi \|_{L^p(\Omega;\R)} := \big ( \E [ |\xi |^p ] \big)^{ 1/p} $ for any $ p \geq 1 $. Let $ x \vee y :=\max \{x, y\} $
and $ x \wedge y :=\min \{x, y\} $
for any $ x, y \in \mathbb{R} $. For notational simplicity,
the letter $ C $ is used to
denote a generic positive constant, which is independent of the
time stepsize and may vary for each appearance.
%
%
Additionally, we suppose that for the coefficient function
$ \varphi : \mathbb{R}  \rightarrow \mathbb{R} $ in \eqref{eq:Ait-Sahalia-model},
there are constants $ M>0 $ and $ \varepsilon_0 > 0 $, such that
\begin{equation}\label{Ass:jump-coefficient-Lipschitz}
|\varphi(x) - \varphi(y) | \leq M |x-y|, \quad \forall x,y>0,
\end{equation}
and
\begin{equation}\label{Ass:jump-coefficient-lower-bound}
x+\varphi(x)>\varepsilon_0 \min\{1,x\}, \quad \forall x>0.
\end{equation}
We comment that the inequality \eqref{Ass:jump-coefficient-lower-bound} naturally
holds when $ x+\varphi(x)>\varepsilon_0 x $ or
$ x+\varphi(x)>\varepsilon_0  $ for some $ \varepsilon_0 > 0 $.
Also, the inequality \eqref{Ass:jump-coefficient-Lipschitz} immediately implies
\begin{equation}
\label{eq:jump-diffusion-growth}
|\varphi(x)| \leq C(1+ | x | ), \quad \forall x > 0,
\end{equation}
where $ C $ is a generic positive constant as explained above.

In the above setting, we attempt to show the existence and uniqueness of
a positive global solution to the problem \eqref{eq:Ait-Sahalia-model}.
The well-posedness of the  Ait-Salalia model without jumps, i.e., $ \varphi \equiv 0 $, has been already 
asserted by \cite[Theorem 2.1]{Szpruch2011Numerical}.
When the jump coefficient $ \varphi $ is a linear function, the well-posedness of the solution to the
model \eqref{eq:Ait-Sahalia-model} has been obtained in \cite[Theorem 2.3]{Deng2018Generalized}.
In what follows, we extend the analysis in \cite{Szpruch2011Numerical,Deng2018Generalized} to 
cope with a more general nonlinear jump coefficient.
%
\begin{proposition} \label{prop:Ait-Sahalia-solution-well-posedness}
	Let conditions \eqref{Ass:jump-coefficient-Lipschitz}, \eqref{Ass:jump-coefficient-lower-bound} hold and let
	the initial data $ X_0 = x_0 >0 $.
	For constants $ a_{-1},\, a_{0}, \, a_{1}, \, a_{2}, \, b >0 $ 
	and $ \gamma, \, \theta > 1 $,
	the problem \eqref{eq:Ait-Sahalia-model}
	admits a unique positive global solution, which 
	almost surely satisfies
	\begin{equation} \label{eq:solution-Ait-Sahalia-model}
	X_{t} = X_{0} + \int_{0}^{t}
	(a_{-1} X_s^{-1} - a_{0}+ a_{1} X_s - a_{2}
	X_s^{\gamma})  \, \dd s
	+\int_{0}^{t} b X_s^{\theta}  \dd W_{s}
	+ \int_{0}^{t} \varphi(X_{s^{-}}) \, \dd N_{s},
	\quad
	t \geq 0.
	\end{equation}
\end{proposition}
\emph{Proof of Proposition \ref{prop:Ait-Sahalia-solution-well-posedness}.}
It is easy to see that the drift and diffusion coefficients of \eqref{eq:solution-Ait-Sahalia-model} are locally Lipschitz continuous in $(0, \infty)$
and that the jump coefficient is globally Lipschitz continuous in $\R$. 
Following the standard truncation arguments \cite{Mao2008Stochastic} and noting $ X_0 = x_0 >0 $,
one can show that there is a unique maximal local solution $X_t, \, t \in\left[0, \tau_{e}\right)$, 
where $\tau_{e}$ is the stopping time of the explosion or first zero time.
To confirm the global solution, we need to prove $ \tau_{e} = \infty $ a.s.
For any sufficiently large positive integer $n$, satisfying $1/n<x_0<n$, we define the stopping times
\begin{equation}\label{eq:defn-stopping-times}
\tau_{n} := \inf \left\{t \in\left[0, \tau_{e}\right): X_t \notin(1 / n, n)\right\},
	\end{equation}
where throughout this paper we set 
 $ \inf ( \emptyset ) = \infty $.
Obviously $ \tau_n $ is increasing as $ n \rightarrow \infty $
and we set $\tau_{\infty} := \lim_{n \rightarrow \infty} \tau_{n}$. 
In view of \eqref{eq:defn-stopping-times}, one knows $ \tau_{\infty} \leq \tau_{e} $ a.s.
If we can prove $\tau_{n} \rightarrow \infty$ a.s.
as $n \rightarrow \infty$, 
then $\tau_{e}=\infty$ a.s. and $x(t) > 0$ a.s. for all $t \geq 0$,
which completes the proof.
To prove $\tau_{\infty}=\infty$ a.s., it suffices to show that 
$\P\left\{\tau_{n} \leq T\right\} \rightarrow 0$ 
as $n \rightarrow \infty$ 
for any constant $T>0,$ which immediately implies $ \P \left\{\tau_{\infty}=\infty\right\}=1$, as required.\\
\indent 
Given a fixed constant  $ \alpha \in (0,1)$,
let us define a function $V \in C^{2}\big((0,\infty), (0,\infty)\big)$ by
\begin{equation} \label{eq:defn-V}
V(x)=x^{\alpha}-\alpha \log x.
\end{equation}
It is easy to check that $V(x) \rightarrow \infty$ 
as $ x \rightarrow \infty$ or $x \rightarrow 0 $ and that
\begin{equation}
\label{eq:VxVxx-computation}
V_{x}(x)=\alpha\big(x^{\alpha-1}-x^{-1}\big),
\quad
V_{xx}(x)=\alpha(\alpha-1) x^{\alpha-2}+\alpha x^{-2}.
\end{equation}
For brevity we write
\begin{equation}\label{eq:ass-drift-diffusion}
\mu(x):=a_{-1} x^{-1}-a_{0}+a_{1} x-a_{2} x^{\gamma}, \quad
\phi(x):= b x^{\theta},
\end{equation}
and introduce the diffusion operator $\mathbb{L} V:(0, \infty) \rightarrow \mathbb{R}$ defined by
\begin{equation} 
\mathbb{L} V(x)=V_{x}(x) \mu(x)+\tfrac{1}{2} V_{xx}(x) \phi^{2}(x).
\end{equation}
Bearing \eqref{eq:VxVxx-computation} and \eqref{eq:ass-drift-diffusion} in mind, one can compute that
\begin{equation} \label{eq:LV-computation}
\begin{split}
\mathbb{L} V(x)  
&  =
\alpha\left(x^{\alpha-1}-x^{-1}\right)\left(a_{-1} x^{-1}-a_{0}+a_{1} x-a_{2} x^{\gamma}\right)+\frac{b^{2}}{2}\left[\alpha(\alpha-1) x^{\alpha-2}+\alpha x^{-2}\right] x^{2 \theta} 
\\
&
=a_{-1} \alpha x^{\alpha-2}-a_{0} \alpha x^{\alpha-1}+a_{1} \alpha x^{\alpha}-a_{2} \alpha x^{\alpha+\gamma-1}-a_{-1} \alpha x^{-2}+a_{0} \alpha x^{-1} 
\\
&
\quad -a_{1} \alpha+a_{2} \alpha x^{\gamma-1}-\frac{b^{2} \alpha(1-\alpha)}{2} x^{\alpha+2 \theta-2}+\frac{b^{2} \alpha}{2} x^{2 \theta-2} 
.
\end{split}
\end{equation}
Using \eqref{Ass:jump-coefficient-lower-bound} and \eqref{eq:jump-diffusion-growth}, one can find a constant $c_{\lambda}$ depending on $\lambda$ such that
\begin{equation} \label{eq: jump-V-computation}
\begin{split}
\lambda(V(x+\varphi(x))-V(x)) 
& =
\lambda[(x+\varphi(x))^{\alpha} - \alpha \log (x+\varphi(x))
-(x^{\alpha} - \alpha \log x)]
\\
&
=
\lambda \Big[ (x+\varphi(x))^{\alpha} -x^{\alpha} - \alpha \log \Big( \frac{x+\varphi(x)}{x} \Big) \Big]
\\
&
\leq
c_{\lambda} ( 1 + x^{\alpha} ),
\quad
\forall x >0.
\end{split}
\end{equation}
Taking \eqref{eq:LV-computation} and  \eqref{eq: jump-V-computation} into account and recalling 
$0<\alpha<1, \gamma>1$, $\theta>1$, $\lambda > 0$,
one can find a constant $ K_1 > 0$ such that 
\begin{equation}
\sup_{x \in(0, \infty)}
\Big(
\mathbb{L} V(x)+\lambda(V(x+\varphi(x))-V(x))
\Big)
 \leq K_1
 < \infty.
\end{equation}
By the It\^{o} formula \cite{Gardon2004approximations} applied to $ V ( x ( t \wedge \tau_{n}) ) $, $ t \in [ 0, T]$, we infer
\begin{equation}
\mathbb{E} 
[
V ( x ( T \wedge \tau_{n}) )
]
 \leq V(x_{0})+ K_1 T
 < \infty,
 \quad
 \forall \, T >0.
\end{equation}
Owing to the definitions \eqref{eq:defn-stopping-times}, \eqref{eq:defn-V}, we deduce from the above estimate that
\begin{equation}
\mathbb{P}\left(\tau_{n} \leq T\right)[V(1/n) \wedge V(n)] 
\leq \mathbb{E} V(x(T \wedge \tau_{n})) 
\leq V\left(x_{0}\right)+ K_1 T < \infty.
\end{equation}
This implies that 
$ \lim_{n \rightarrow \infty}
\mathbb{P}\left(\tau_{n} \leq T\right) = 0$ for any constant $T>0$ and the proof is thus complete.
\qed
In the following error analysis, the moment bounds of the analytic solution are frequently used.
The next lemma indicates when $ p $-th moments of the solution to \eqref{eq:Ait-Sahalia-model} are bounded.
\begin{lemma} \label{lem:solution-moment-bound}
	Let all conditions in Proposition \ref{prop:Ait-Sahalia-solution-well-posedness} hold and let $ \{ X_t \}_{ t \geq 0 } $ 
	be the unique solution to \eqref{eq:Ait-Sahalia-model}, given by \eqref{eq:solution-Ait-Sahalia-model}.
	If one of the
	following two conditions holds: \\
	(i) $ p \geq 2 $ when $ \gamma+1> 2\theta $; \\
	(ii) $ 2 \leq p < \frac{2 a_2 + b^2}{b^2} $ when $ \gamma+1=2\theta $,\\
	then
	\begin{equation}\label{eq:p-moment-bound}
	\sup _{t \in[0, \infty)} \E [ |X_{t} |^{p} ]<\infty.
	\end{equation}
\end{lemma}
\emph{Proof of Lemma \ref{lem:solution-moment-bound}.}
For a sufficiently large positive integer
$ n $ satisfying $ \frac{1}{n}< x_0 < n $,
we define the stopping time
\begin{equation}
\tau_{n} := \inf  \{t \in [0, \infty ) :
X_t \notin(1 / n, n) \}.
\end{equation}
Also, we define
$ V_1 $ : $ \mathbb{R}_{+} \times [0,\infty) \rightarrow \mathbb{R}_{+} $
as follows
\begin{equation}
V_1(x,t) : = e^{t} x^{p}, \quad x \in \mathbb{R}_{+}, \, t \in [0,\infty).
\end{equation}
We compute that
\begin{equation}\label{eq:operator-L-computation}
\begin{split}
& \quad
\mathbb{L} V_1(x,t) + \lambda (V_1
( x + \varphi(x),t) - V_1(x,t)) \\
& = e^{t} \big[ p x^{p-1}
(a_{-1} x^{-1} - a_{0} + a_{1} x - a_{2}
x^{\gamma}) + \tfrac{1}{2} b^{2} p(p-1)
x^{p-2+2 \theta} \big] \\
& \quad + \lambda e^{t}((x + \varphi(x))^{p}
- x^{p} ) \\
&=
e^{t} \big[a_{-1} p x^{p-2} - a_{0}p x^{p-1}
+ a_{1} p x^{p} - a_{2}p x^{p + \gamma -1}
+ \tfrac{1}{2} b^2 p(p-1)x^{p-2+2\theta} \\
& \quad
+  \lambda ((x + \varphi(x))^{p} - x^{p} )\big],
\end{split}
\end{equation}
where  the diffusion operator
$ \mathbb{L} V_1 : \mathbb{R}_{+} \times [0, \infty)
\rightarrow \mathbb{R}_{+} $
is defined by
$$
\mathbb{L} V_1(x,t) := \tfrac{\partial V_1(x,t)}{\partial x}
\mu(x)
+
\tfrac{1}{2} \tfrac{\partial^2 V_1(x,t)}{\partial x^2} \phi^2 (x).
$$
In light of condition (i) or (ii) as well as \eqref{Ass:jump-coefficient-lower-bound}, \eqref{eq:jump-diffusion-growth},
one can find a constant $ K_2 > 0$, such that
\begin{equation}\label{ineq:operator-L-bound}
\mathbb{L} V_1(X_t) + \lambda (V_1
( X_t + \varphi(X_t)) - V_1(X_t)) \leq K_2 e^t.
\end{equation}
Indeed, in the case (i),
we can directly get $ p + \gamma -1 > p - 2 + 2 \theta $.
Thus, it is easy to see that the highest power of $ x $ is
$ p + \gamma -1 $ in \eqref{eq:operator-L-computation}. As a result of $ a_2 >0 $,
there is a constant $ K_2 >0 $ such that \eqref{ineq:operator-L-bound} is fulfilled.
In the case (ii),
one can derive that $ p + \gamma -1 = p - 2 + 2 \theta $.
Furthermore, the inequality $ \tfrac{1}{2} b^2 p (p-1) - a_2 p < 0 $
holds under the condition  $ p < \tfrac{2a_2+b^2}{b^2} $.
Similar to the above analysis, we can arrive at the same conclusion.
By the It\^{o} formula \cite{Gardon2004approximations}, for any $ t \geq 0 $,
\begin{equation}
\E[e^{t \wedge \tau_{n}} X_{t \wedge \tau_{n}}^{p}
]
\leq x_{0}^{p} + K_2 e^{t}.
\end{equation}
Letting $n  \rightarrow \infty$ and applying Fatou's lemma,
we obtain
\begin{equation}
\E [|X_t|^{p}] \leq \frac{x_0^{p}}{e^{t}} + K_2, \quad t \geq 0.
\end{equation}
The proof of Lemma \ref{lem:solution-moment-bound}
is thus completed.
\qed
The next lemma gives the inverse moment bounds of the solution to \eqref{eq:Ait-Sahalia-model}, which is also required in the error analysis.
\begin{lemma}\label{lem:solution-inverse-moment-bound}
	Let all conditions in Proposition \ref{prop:Ait-Sahalia-solution-well-posedness} hold with $ \gamma +1 \geq 2 \theta $ and let $ \{ X_t \}_{ t \geq 0 } $ 
	be the unique solution to \eqref{eq:Ait-Sahalia-model}, given by \eqref{eq:solution-Ait-Sahalia-model}.
	Then for any $ p \geq \max\{ 1, \gamma - 1\} $
	it holds
	\begin{equation}\label{eq:inverse-p-moment-bound}
	\sup _{t \in[0, \infty)} \E [ |X_{t} |^{-p} ]<\infty.
	\end{equation}
\end{lemma}
\emph{Proof of Lemma \ref{lem:solution-inverse-moment-bound}.}
Define
$ V_2 $ : $ \mathbb{R}_{+} \times [0,\infty) \rightarrow \mathbb{R}_{+} $
as follows
\begin{equation}
V_2 ( x, t) := e^{t} x^{-p}, \quad x \in \mathbb{R}_{+}, \,\,\, t \in [0,\infty).
\end{equation}
Here, $ \tau_{n} $ and the functional $ \mathbb{L} $
have been defined in the proof of Lemma \ref{lem:solution-moment-bound}.
Then we have
\begin{equation}\label{eq:operator-L-computationoperator-inverde-moment}
\begin{split}
& \quad
\mathbb{L} V_2(x,t) + \lambda (V_2
( x + \varphi(x),t) - V_2(x,t)) \\
&=
e^{t} \big[- a_{-1} p x^{-p-2} + a_{0}p x^{-p-1}
- a_{1} p x^{-p} + a_{2}p x^{-p + \gamma -1}
+ \tfrac{1}{2} b^2 p(p+1)x^{-p-2+2\theta} \\
& \quad
+  \lambda ((x + \varphi(x))^{-p} - x^{-p} )\big].
\end{split}
\end{equation}
Taking $ \gamma+1  \geq 2\theta  $ and $ \gamma \leq p+1 $
into account promises that the highest power of
$ x $ on the right-hand side of \eqref{eq:operator-L-computationoperator-inverde-moment} is $ -p + \gamma -1 \leq 0 $.
Also, it is not difficult to check that $ -p-2  $ is the lowest power of $ x $
in \eqref{eq:operator-L-computationoperator-inverde-moment}. Due to the negative coefficient of $ x^{-p-2} $, i.e, $ -a_{-1}p < 0 $,
there exists a constant $ K_3 > 0 $ such that
\begin{equation}\label{ineq:operator-L-bound-inverse-moment}
\mathbb{L} V_2(X_t) + \lambda \big(V_2
( X_t + \varphi(X_t)) - V_2(X_t)\big) \leq K_3 e^t.
\end{equation}
The remaining proof is similar to that of Lemma
\ref{lem:solution-moment-bound}
and thus omitted.
\qed

\section{Mean-square convergence rate of the backward Euler method for the generalized Ait-Sahalia-type
	rate model with Poisson jumps}
In this section we aim to propose and analyze the backward Euler method (BEM) for strong approximations of the Ait-Sahalia model 
with Poisson jumps  \eqref{eq:Ait-Sahalia-model}. 
Given $T \in (0, \infty)$ and $N \in \N$,  we construct a uniform mesh on the interval $[0, T]$ with the uniform stepsize $ h = \tfrac{T}{N}$.
Based on the uniform mesh, we propose a numerical scheme for \eqref{eq:Ait-Sahalia-model} as follows:
\begin{equation}\label{eq:BEM-numerical-approximation}
Y_{n}= Y_{n-1} +  h [a_{-1} Y_{n}^{-1}-a_{0}
+ a_{1} Y_{n}-a_{2} Y_{n}^{\gamma} ]
+ b Y_{n-1}^{\theta} \Delta W_{n-1} + \varphi(Y_{n-1})\Delta N_{n-1},
\quad Y_{0}=X_{0},
\end{equation}
where $ \Delta W_{n-1} := W_{t_n} - W_{t_{n-1}} $,
$ \Delta N_{n-1} := N_{t_n} - N_{t_{n-1}} $, $ n \in \{ 1, 2, ..., N \}, N \in \N$.
Next we shall check that the numerical approximations produced by \eqref{eq:BEM-numerical-approximation} 
are well defined and preserves positivity of the original model \eqref{eq:Ait-Sahalia-model}.
\begin{lemma}\label{lem:numerical-solution-well-posedness}
	Let all conditions in Proposition \ref{prop:Ait-Sahalia-solution-well-posedness} hold.
	For any $ h \leq \frac{1}{ a_{1}} $, the BEM \eqref{eq:BEM-numerical-approximation} is well-defined in the sense that it admits a
	unique positive solution.
\end{lemma}
\emph{Proof of Lemma \ref{lem:numerical-solution-well-posedness}.}
The proof is similar to that of \cite[Lemma 3.1]{Szpruch2011Numerical}.
\qed
%
To carry out the error analysis of the BEM for the model, we treat the standard
case $ \kappa + 1 > 2 \rho $ and the critical case $ \kappa + 1 = 2 \rho $ separately.
\subsection{The case $\gamma+1>2 \theta$}
This subsection attempts to recover the expected convergence rate of the BEM for  
the Ait-Sahalia model in the case $ \gamma+1>2 \theta $. Before
proceeding further, we present an important lemma,
which plays an essential role in deriving the mean-square convergence rate of BEM.
\begin{lemma}\label{lem:solution-estimation}
Let all conditions in Proposition \ref{prop:Ait-Sahalia-solution-well-posedness} hold with 
$ \gamma + 1 > 2 \theta $ and let $ \{ X_t \}_{ t \geq 0 } $ 
	be the process defined by \eqref{eq:solution-Ait-Sahalia-model}.
	Then it holds that
	\begin{equation}\label{eq:solution-esti}
	\|X_{t}-X_{s} \|_{L^{2} (\Omega ; \mathbb{R})}
	\leq C|t-s|^{\frac{1}{2}}, \quad 0 \leq t,s \leq T.
	\end{equation}
\end{lemma}
\emph{Proof of Lemma \ref{lem:solution-estimation}.}
Recall that the coefficients of \eqref{eq:Ait-Sahalia-model} are defined as follows:
\begin{equation}
\mu(x):=a_{-1} x^{-1}-a_{0}+a_{1} x-a_{2} x^{\gamma}, \quad
\phi(x):= b x^{\theta}.
\end{equation}
Letting $ t > s $ without loss of generality,
then one can easily check that
\begin{equation}
\E [ | X_{t}-X_{s}|^2]
= \E \Big[ \Big| \int_{s}^{t} \mu(X_r) \, \dd r +\int_{s}^{t}
\phi( X_{r}) \, \dd W_{r}
+ \int_{s}^{t} \varphi(X_{r^{-}})  \, \dd \tilde{N}_r
+ \lambda \int_{s}^{t} \varphi(X_{r^{-}}) \, \dd r \Big|^2 \Big].
\end{equation}
Employing the Young inequality, the H\"older inequality,
the It\^{o} isometry and Lemmas
\ref{lem:solution-moment-bound},
\ref{lem:solution-inverse-moment-bound},
one can arrive at
\begin{equation}\label{eq:solution-order}
\begin{split}
\E [ | X_{t}-X_{s}|^2]
\leq &
4 |t-s| \int_{s}^{t} \E [| \mu(X_r) |^2 ] \,\dd r
+ 4 \int_{s}^{t}
\E[| \phi(X_{r})|^2 ] \, \dd r \\
& +
4 \lambda \int_{s}^{t} \E [ |\varphi(X_{r^{-}})|^2 ] \, \dd r
+
4 \lambda^2 |t-s|
\int_{s}^{t} \E[| \varphi(X_{r^{-}}) |^2 ] \, \dd r \\
\leq &
C | t - s |^2 + C |t-s|  \\
\leq &
C | t - s |,
\end{split}
\end{equation}
as required.
\qed
Equipped with the above lemmas, we are ready to derive
the mean-square convergence rate of order one-half for the scheme.
\begin{theorem}\label{thm:convergence-rate}
Let all conditions in Proposition \ref{prop:Ait-Sahalia-solution-well-posedness} hold with $\gamma + 1 > 2 \theta$.
	Let $\{ X_{t} \}_{0 \leq t \leq T }$
	and $\{ Y_{n} \}_{0 \leq n \leq N } $ be solutions to \eqref{eq:Ait-Sahalia-model}
	and \eqref{eq:BEM-numerical-approximation}, respectively. 
	Then for any $h \in (0,  \frac{1}{2L} )$ it holds that
	\begin{equation}\label{eq:convergence-rate}
	\sup _{N \in \mathbb{N}} \sup _{0 \leq n \leq N}
	\|Y_{n} -  X_{t_{n}} \|_{L^{2}
		(\Omega ; \mathbb{R})} \leq C h^{\frac{1}{2 }},
	\end{equation}
	where $ L := a_1 + \frac{(q-1)b^2 \theta^2(\gamma + 1 - 2 \theta)}
	{2(\gamma-1)}
	\big( \frac{(q-1)b^2 \theta^2(\theta-1)}{ a_2 \gamma(\gamma-1)}
	\big)^{\frac{2\theta-2}{\gamma + 1 - 2 \theta}}$ for any $ q > 2 $.
\end{theorem}
\emph{Proof of Theorem \ref{thm:convergence-rate}.}
By the definition \eqref{eq:Ait-Sahalia-model}, we learn that
\begin{equation}\label{eq:one-step-solution}
X_{t_n}= X_{t_{n-1}}+ \int^{t_{n}}_{t_{n-1}}
(a_{-1} X_s^{-1} - a_{0}+ a_{1} X_s - a_{2} X_s^{\gamma} ) \, \dd s
+ \int^{t_{n}}_{t_{n-1}}b X_s^{\theta} \dd W_s
+ \int^{t_{n}}_{t_{n-1}} \varphi(X_{s^{-}}) \, \dd N_s.
\end{equation}
Subtracting \eqref{eq:one-step-solution} from \eqref{eq:BEM-numerical-approximation}
leads to
\begin{equation}\label{eq:solution-numerical-difference}
\begin{split}
E_{n} = &E_{n-1} +
h \Delta \mu_{n}
+ \Delta \phi_{n-1} \Delta W_{n-1}
+ \Delta \varphi_{n-1} \Delta N_{n-1}
+ \mathcal{M}_{t_n} ,
\end{split}
\end{equation}
where for short we denote
\begin{equation}
\begin{split}
& E_{n} := Y_{n} - X_{t_n}, \quad
\Delta \mu_{n} := \mu(Y_n) - \mu(X_{t_n}), \\
& \Delta \phi_{n-1} := \phi(Y_{n-1}) - \phi(X_{t_{n-1}}), \quad
\Delta \varphi_{n-1} := \varphi(Y_{n-1}) - \varphi(X_{t_{n-1}}),\\
& \mathcal{M}_{t_n} := \int_{t_{n-1}}^{t_{n}} \mu (X_{t_{n}} ) - \mu (X_{s} ) \, \dd s
+  \int_{t_{n-1}}^{t_{n}} \phi(X_{t_{n-1}}) - \phi(X_{s})  \, \dd W_{s}
+  \int_{t_{n-1}}^{t_{n}} \varphi(X_{t_{n-1}}) - \varphi(X_{s^{-}}) \, \dd N_{s}.
\end{split}
\end{equation}
Further, one can deduce from \eqref{eq:solution-numerical-difference} that
\begin{equation}\label{eq:En-hDeltamu_n}
| E_{n} - h \Delta \mu_n |^2
= | E_{n-1} + \Delta \phi_{n-1} \Delta W_{n-1}
+ \Delta \varphi_{n-1} \Delta N_{n-1}
+ \mathcal{M}_{t_{n}}|^2.
\end{equation}
By direct calculation, one can recast \eqref{eq:En-hDeltamu_n} as
\begin{equation}\label{eq:En-expansion}
\begin{split}
& | E_{n} |^2 + |h \Delta \mu_n |^2
\\
= &
| E_{n-1} |^2 + | \Delta \phi_{n-1} \Delta W_{n-1}|^2
+ | \Delta \varphi_{n-1} \Delta N_{n-1}|^2
+ | \mathcal{M}_{t_{n}} |^2
+ 2 \langle E_{n-1} , \Delta \phi_{n-1} \Delta W_{n-1} \rangle \\
&+2 \langle E_{n-1} , \Delta \varphi_{n-1} \Delta N_{n-1} \rangle
+ 2 \langle E_{n-1} , \mathcal{M}_{t_{n}} \rangle
+ 2 \langle \Delta \phi_{n-1} \Delta W_{n-1} ,
\Delta \varphi_{n-1} \Delta N_{n-1} \rangle \\
&+ 2 \langle \Delta \phi_{n-1} \Delta W_{n-1} , \mathcal{M}_{t_{n}} \rangle
+ 2 \langle \Delta \varphi_{n-1} \Delta N_{n-1} , \mathcal{M}_{t_{n}} \rangle
+ 2h \langle E_{n} , \Delta \mu_n \rangle.
\end{split}
\end{equation}
Taking expectations on both sides of \eqref{eq:En-expansion},
one can easily arrive at
\begin{equation}\label{eq:En-expansion-expectation}
\begin{split}
& \E [| E_{n} |^2] + h^2 \E[| \Delta \mu_n |^2] \\
= &
\E[| E_{n-1} |^2] + \E[| \Delta \phi_{n-1} \Delta W_{n-1}|^2]
+ \E[| \Delta \varphi_{n-1} \Delta N_{n-1}|^2]
+ \E [| \mathcal{M}_{t_{n}} |^2] \\
& + 2 \E[ \langle E_{n-1} , \Delta \phi_{n-1} \Delta W_{n-1} \rangle ]
+2 \E[\langle E_{n-1} , \Delta \varphi_{n-1} \Delta N_{n-1} \rangle] \\
& + 2 \E[\langle E_{n-1} , \mathcal{M}_{t_{n}} \rangle]
+ 2 \E[\langle \Delta \phi_{n-1} \Delta W_{n-1} ,
\Delta \varphi_{n-1} \Delta N_{n-1} \rangle] \\
& + 2 \E [\langle \Delta \phi_{n-1} \Delta W_{n-1} ,
\mathcal{M}_{t_{n}} \rangle]
+ 2 \E[\langle \Delta \varphi_{n-1} \Delta N_{n-1} ,
\mathcal{M}_{t_{n}} \rangle]
+ 2h \E[\langle E_{n} , \Delta \mu_n \rangle].
\end{split}
\end{equation}
Before proceeding further, we claim that
\begin{equation}\label{eq:remainder-term-bound}
\E [| \mathcal{M}_{t_{n}} |^2] < \infty .
\end{equation}
Indeed, using the Young inequality, the H\"older inequality, the It\^{o} isometry
and Lemmas \ref{lem:solution-moment-bound}, \ref{lem:solution-inverse-moment-bound} shows
\begin{equation}\label{eq:proof-remainder-term-bound}
\begin{split}
\E [| \mathcal{M}_{t_{n}} |^2]
=&
\E \Big[ \Big| \int_{t_{n-1}}^{t_n}
\mu(X_{t_n}) - \mu(X_s) \,\dd s
+ \int_{t_{n-1}}^{t_n}
\phi(X_{t_{n-1}}) - \phi(X_s)\, \dd W_s \\
&  + \int_{t_{n-1}}^{t_n}
\varphi(X_{t_{n-1}}) - \varphi(X_{s^{-}}) \, \dd N_s \Big|^2 \Big] \\
\leq&
4 \E \Big[ \Big| \int_{t_{n-1}}^{t_n}
\mu(X_{t_n}) - \mu(X_s) \,\dd s \Big|^2 \Big]
+ 4 \E \Big[ \Big| \int_{t_{n-1}}^{t_n}
\phi(X_{t_{n-1}}) - \phi(X_s) \, \dd W_s \Big|^2 \Big]\\
& +
4 \E \Big[ \Big| \int_{t_{n-1}}^{t_n}
\varphi(X_{t_{n-1}}) - \varphi(X_{s^{-}}) \, \dd \tilde{N}_s \Big|^2 \Big]
+ 4 \E \Big[ \lambda^2 \Big| \int_{t_{n-1}}^{t_n}
\varphi(X_{t_{n-1}}) - \varphi(X_{s^{-}}) \, \dd s \Big|^2 \Big] \\
\leq &
4 h \int_{t_{n-1}}^{t_n}
\E \big[| \mu(X_{t_n}) - \mu(X_s)|^2 \big]  \dd s
+ 4    \int_{t_{n-1}}^{t_n}
\E \big[ | \phi(X_{t_{n-1}}) - \phi(X_s)|^2 \big] \dd s\\
&+
4 \lambda \int_{t_{n-1}}^{t_n}
\E \big[ | \varphi(X_{t_{n-1}}) - \varphi(X_{s^{-}})|^2 \big] \dd s
+ 4 \lambda^2 h \int_{t_{n-1}}^{t_n}
\E \big[ | \varphi(X_{t_{n-1}}) - \varphi(X_{s^{-}})|^2 \big] \dd s  \\
< & \infty.
\end{split}
\end{equation}
In addition, we note that for any
$ \gamma +1 > 2 \theta $ and $ q> 2 $,
\begin{equation}
\sup_{x>0}\Big(\mu^{\prime}(x) + \tfrac{q-1}{2}(\phi'(x))^{2}\Big)
=\sup_{x>0}\Big(-a_{-1} x^{-2}+ a_{1}- a_{2} \gamma x^{\gamma-1}
+ \tfrac{q-1}{2} b^{2} \theta^{2}
x^{2 \theta-2}\Big) < \infty, 
\end{equation}
which in turn implies
\begin{equation}\label{eq:couple-monotonicity-condition}
\langle x-y, \mu(x)- \mu(y)\rangle + \tfrac{q-1}{2} |\phi(x)-\phi(y)|^{2}
\leq L |x-y|^{2}, \quad \forall x, y \in {\R}_{+},
\end{equation}
where $ L :=
a_1 + \sup_{x \in (0,\infty)} \big(\tfrac{q-1}{2}b^2\theta^2 x^{2\theta-2}
-\gamma a_2 x^{\gamma-1}\big)
= a_1 + \frac{(q-1)b^2 \theta^2(\gamma + 1 - 2 \theta)}
{2(\gamma-1)}
\big( \frac{(q-1)b^2 \theta^2(\theta-1)}{ a_2 \gamma(\gamma-1)}
\big)^{\frac{2\theta-2}{\gamma + 1 - 2 \theta}} $.
As a direct result of \eqref{eq:couple-monotonicity-condition},
one can infer from \eqref{eq:En-expansion-expectation} that
\begin{equation}\label{eq:En-expectation-expansion-esti}
\begin{aligned}
& (1 - 2hL) \E [| E_{n} |^2] + h^2 \E[| \Delta \mu_n |^2] \\
\leq &
\E[| E_{n-1} |^2] + \E[| \Delta \phi_{n-1} \Delta W_{n-1}|^2]
+ \E[| \Delta \varphi_{n-1} \Delta N_{n-1}|^2]
+ \E [| \mathcal{M}_{t_{n}} |^2] \\
& + 2 \E[ \langle E_{n-1} , \Delta \phi_{n-1} \Delta W_{n-1} \rangle ]
+2 \E[\langle E_{n-1} , \Delta \varphi_{n-1} \Delta N_{n-1} \rangle]
+ 2 \E[\langle E_{n-1} , \mathcal{M}_{t_{n}} \rangle] \\
& + 2 \E[\langle \Delta \phi_{n-1} \Delta W_{n-1} ,
\Delta \varphi_{n-1} \Delta N_{n-1} \rangle] \\
& + 2 \E [\langle \Delta \phi_{n-1} \Delta W_{n-1} ,
\mathcal{M}_{t_{n}} \rangle]
+ 2 \E[\langle \Delta \varphi_{n-1} \Delta N_{n-1} ,
\mathcal{M}_{t_{n}} \rangle].
\end{aligned}
\end{equation}
In the next step, we declare that
\begin{equation}\label{eq:E_k-drift-bound}
\E [ |E_k|^2 ] < \infty, \quad
\E [ |\Delta \mu_k |^2 ] < \infty, \quad
\forall \, k \in \{ 0, 1, 2, ..., N \},
\,
N\in \N,
\end{equation}
whose proof relies on the mathematical induction argument.
In fact, for $ k = 0 $,
$ Y_0 = X_0 $ and thus $ E_{0}
= \Delta \mu_{0} = 0 $,
the assertion naturally holds.
Next we assume, for some $ n \leq N $,
\begin{equation}\label{eq:En-drift-bound}
\E [ | E_{n-1} |^2 ] < \infty,
\quad
\E [ |\Delta \mu_{n-1} |^2 ] < \infty.
\end{equation}
In view of \eqref{Ass:jump-coefficient-Lipschitz} and \eqref{eq:couple-monotonicity-condition},
one can show that
\begin{equation}\label{eq:diffusion-jump-bound}
\begin{split}
&\E[|\Delta \varphi_{n-1}|^{2}]
\leq
M^2 \E[|E_{n-1} |^{2}]< \infty, \\
\quad
&\E[|\Delta \phi_{n-1}|^{2}]
\leq
\tfrac{1}{q-1}(2L+1) \E[ |E_{n-1}|^{2}]
+ \tfrac{1}{q-1} \E[|\Delta \mu_{n-1} |^{2}] < \infty .
\end{split}
\end{equation}
Employing these bounded moments above enables us to
deduce that
\begin{equation}\label{eq:Ito-isometry}
\begin{split}
& \E [ |\Delta \phi_{n-1} \Delta W_{n-1} |^{2} ]
= h \E [ |\Delta \phi_{n-1} |^{2} ], 
\\
&
\E [ \langle E_{n-1},
\Delta \phi_{n-1}
\Delta W_{n-1} \rangle ] = 0 , \\
& \E [\langle \Delta \phi_{n-1} \Delta W_{n-1} ,
\Delta \varphi_{n-1} \Delta N_{n-1} \rangle] = 0,
\end{split}
\end{equation}
and
\begin{equation}\label{eq:Poisson-Ito-isometry}
\begin{split}
\E [ |\Delta \varphi_{n-1} \Delta N_{n-1} |^{2} ]
&= \E [ | \Delta \varphi_{n-1} |^{2} ]
\cdot \E [ | \Delta N_{n-1}  |^{2} ] \\
& = \E [ | \Delta \varphi_{n-1} |^{2} ]
\cdot \big( \E [ | \Delta \tilde{N}_{n-1}  |^{2} ]
+ \E [ | \lambda h |^2 ]
+ 2 \E [ \lambda h \Delta \tilde{N}_{n-1}] \big)  \\
& = ( \lambda h + \lambda^2 h^2) \E [ | \Delta \varphi_{n-1} |^{2} ],
\end{split}
\end{equation}
\begin{equation}\label{eq:property-compensated-Poisson-process}
\begin{split}
\E [ \langle E_{n-1} , \Delta \varphi_{n-1} \Delta N_{n-1} \rangle ]
&= \E [ \langle E_{n-1} , \Delta \varphi_{n-1} \Delta \tilde{N}_{n-1} \rangle ]
+ \lambda h \E [ \langle E_{n-1} , \Delta \varphi_{n-1}  \rangle ] \\
&= \lambda h \E [ \langle E_{n-1} , \Delta \varphi_{n-1}  \rangle ],
\end{split}
\end{equation}
where $ \Delta \tilde{N}_{n-1} :=  \tilde{N}_{t_{n}} - \tilde{N}_{t_{n-1}}$.
Therefore, using \eqref{eq:proof-remainder-term-bound},
\eqref{eq:En-expectation-expansion-esti}, \eqref{eq:En-drift-bound}, \eqref{eq:diffusion-jump-bound}
and the Cauchy-Schwarz inequality,
we have
\begin{equation}
\begin{aligned}
& (1 - 2hL) \E [| E_{n} |^2] + h^2 \E[| \Delta \mu_n |^2] \\
\leq&
\E[ |E_{n-1}|^2 ] + h \E[ |\Delta \phi_{n-1} |^2]
+ (\lambda h + \lambda^2 h^2) \E [ |\Delta \varphi_{n-1} |^2 ]
+ \E [ | \mathcal{M}_{t_n} |^2 ] \\
& + 2 \lambda h \E [ \langle E_{n-1} , \Delta \varphi_{n-1} \rangle ]
+  \E [ | E_{n-1}|^2] + \E[ |\mathcal{M}_{t_n}|^2 ]
+ h \E [ | \Delta \phi_{n-1} |^2 ]  \\
& + \E [ | \mathcal{M}_{t_n} |^2 ]
+ (\lambda h + \lambda^2 h^2) \E [ | \Delta \varphi_{n-1} |^2 ]
+ \E [ | \mathcal{M}_{t_n} |^2 ]  \\
\leq &
(2+ \lambda h) \E[| E_{n-1} |^2] + 2 h \E[| \Delta \phi_{n-1} |^2]
+ (3\lambda h + 2 \lambda^2 h^2 )\E[| \Delta \varphi_{n-1} |^2]
+ 4 \E [| \mathcal{M}_{t_{n}} |^2] < \infty.
\end{aligned}
\end{equation}
This implies
\begin{equation}
\E [| E_{n} |^2] < \infty ,  \quad
\E[| \Delta \mu_n |^2]  < \infty,
\end{equation}
and the claim \eqref{eq:E_k-drift-bound} is thus validated based on the mathematical induction argument.
Similar to \eqref{eq:diffusion-jump-bound}, we use \eqref{eq:E_k-drift-bound} to infer
\begin{equation}\label{eq:diffusion-jump-bounds}
\begin{split}
\E[|\Delta \varphi_{k}|^{2}]
< \infty, 
\quad
\E[|\Delta \phi_{k} |^{2} ]
< \infty,
\quad
\forall \, k \in \{ 0, 1, 2, ..., N \},
\,
N\in \N.
\end{split}
\end{equation}
At the moment we are well-prepared to prove the desired error estimate \eqref{eq:convergence-rate}.
Using \eqref{Ass:jump-coefficient-Lipschitz}, \eqref{eq:couple-monotonicity-condition}, \eqref{eq:Ito-isometry},
\eqref{eq:Poisson-Ito-isometry}, \eqref{eq:property-compensated-Poisson-process},
the Cauchy-Schwarz inequality and the properties
of the conditional expectation,
we deduce from \eqref{eq:En-expansion-expectation} that
\begin{equation}\label{eq:convergen-error-order}
\begin{split}
& \E [ |E_{n} |^{2} ] - \E [ |E_{n-1} |^{2} ] \\
\leq &
h \E [ | \Delta \phi_{n-1} |^2 ]
+ (\lambda h + \lambda^2 h^2) \E [ | \Delta \varphi_{n-1} |^2 ]
+ \E[| \mathcal{M}_{t_n}  |^{2}]
+ \lambda h \E [ |E_{n-1} |^{2} ] \\
& + \lambda h \E [ | \Delta \varphi_{n-1} |^2 ]
+ 2 \E [\langle E_{n-1} ,
\mathcal{M}_{t_{n}} \rangle ]
+ (q-2)h \E [ | \Delta \phi_{n-1} |^2 ]
+ \tfrac{1}{q-2} \E[| \mathcal{M}_{t_n}  |^{2}] \\
&+ (\lambda h + \lambda^2 h^2) \E [ | \Delta \varphi_{n-1} |^2 ]
+ \E[| \mathcal{M}_{t_n}  |^{2}]
+ 2 h \E [\langle E_{n} ,\Delta \mu_{n} \rangle ]
\\
\leq &
\lambda h \E [ |E_{n-1} |^{2} ]
+ (q-1) h \E [ | \Delta \phi_{n-1} |^2 ]
+ (3 \lambda h + 2\lambda^2 h^2) \E [ | \Delta \varphi_{n-1} |^2 ]\\
&
+ 2 \E [\langle E_{n-1} , \E(\mathcal{M}_{t_{n}}
|\mathcal{F}_{t_{n-1}})\rangle ]
+  \tfrac{2q-3}{q-2} \E[| \mathcal{M}_{t_n}  |^{2}]
+ 2 h \E [\langle E_{n} ,\Delta \mu_{n} \rangle ] \\
\leq &
\big( \lambda h + (3\lambda h + 2\lambda^2 h^2) M^2  \big)
\E [ |E_{n-1} |^{2} ]
+ (q-1) h \E [ | \Delta \phi_{n-1} |^2 ]
+  h \E [ | E_{n-1}|^2 ]  \\
& + h^{-1} \E [ | \E(\mathcal{M}_{t_{n}}
|\mathcal{F}_{t_{n-1}})|^2 ]
+ \tfrac{2q-3}{q-2} \E[| \mathcal{M}_{t_n}  |^{2}]
+ 2 L h \E [ |E_{n} |^{2} ] - (q-1) h \E [ | \Delta \phi_{n} |^2 ] \\
= &
\big( (\lambda + 1) h + (3\lambda h + 2\lambda^2 h^2) M^2  \big)
\E [ |E_{n-1} |^{2} ]
+ (q-1) h \E [ | \Delta \phi_{n-1} |^2 ] \\
& + h^{-1} \E [ | \E(\mathcal{M}_{t_{n}}
|\mathcal{F}_{t_{n-1}})|^2 ]
+ \tfrac{2q-3}{q-2} \E[| \mathcal{M}_{t_n}  |^{2}]
+ 2 L h \E [ |E_{n} |^{2} ] - (q-1) h \E [ | \Delta \phi_{n} |^2 ].
\end{split}
\end{equation}
Summing both sides of the above inequality from
$  1 $ to $ n $, we get
\begin{equation}\label{eq:En-mean-square-expectation-esti}
\begin{split}
\E [ |E_{n} |^{2} ]
\leq &
2Lh\E[ |E_n|^2] + \big( 2Lh +(\lambda + 1) h
+ (3\lambda h + 2\lambda^2 h^2) M^2  \big)
\sum_{k=0}^{n-1}\E [ |E_{k} |^{2} ] \\
& + h^{-1} \sum_{k=1}^{n} \E [ | \E(\mathcal{M}_{t_{k}}
|\mathcal{F}_{t_{k-1}})|^2 ]
+ \tfrac{2q-3}{q-2} \sum_{k=1}^{n} \E[| \mathcal{M}_{t_k}  |^{2}]
- (q-1) h \E [ | \Delta \phi_{n} |^2 ] \\
\leq &
2Lh\E[ |E_n|^2] + C h \sum_{k=0}^{n-1}\E [ |E_{k} |^{2} ]
\\
& + h^{-1} \sum_{k=1}^{n} \E [ | \E(\mathcal{M}_{t_{k}}
|\mathcal{F}_{t_{k-1}})|^2 ]
+ \tfrac{2q-3}{q-2} \sum_{k=1}^{n} \E[| \mathcal{M}_{t_k}  |^{2}].
\end{split}
\end{equation}
Noting $ 1 - 2 Lh > \zeta > 0 $ for some $ \zeta > 0 $ and invoking the Gronwall inequality gives
\begin{equation} \label{eq:error-bound}
\E [ |E_{n} |^{2} ]
\leq
C \Big( h^{-1} \sum_{j=1}^{n} \E[ |\E (\mathcal{M}_{t_{j}}
| \mathcal{F}_{t_{j-1}})|^2  ]
+   \sum_{j=1}^{n} \E[| \mathcal{M}_{t_j} |^{2}] \Big) .
\end{equation}
Accordingly,
it remains to estimate
$ \E [ |\mathcal{M}_{t_{j}} |^{2} ] $ and
$ \E [ |\E (\mathcal{M}_{t_{j}} | \mathcal{F}_{t_{j-1}} ) |^{2} ] $
before attaining the mean-square convergence rate.
An elementary inequality gives
\begin{equation}\label{eq:remainder-term-esti}
\begin{aligned}
\|\mathcal{M}_{t_{j}} \|_{L^{2} (\Omega ; \mathbb{R} )}
= &
\Big\| \int_{t_{j-1}}^{t_{j}} \mu (X_{t_{j}} )
- \mu (X_{s} )\, \dd s
+  \int_{t_{j-1}}^{t_{j}} \phi(X_{t_{j-1}}) - \phi(X_{s}) \, \dd W_{s} \\
& +
\int_{t_{j-1}}^{t_{j}} \varphi(X_{t_{j-1}}) - \varphi(X_{s^{-}}) \, \dd N_{s}
\Big\|_{L^{2} (\Omega ; \mathbb{R} )} \\
\leq &
C \int_{t_{j-1}}^{t_{j}}\|X_{t_{j}}^{-1}
- X_{s}^{-1}\|_{L^{2}(\Omega ; \mathbb{R})}
+ \|X_{t_{j}}- X_{s}\|_{L^{2}(\Omega ;\mathbb{R})}
+ \|X_{t_{j}}^{\gamma}-X_{s}^{\gamma}
\|_{L^{2}(\Omega ; \mathbb{R})} \, \dd s \\
& +
b \Big\| \int_{t_{j-1}}^{t_{j}} X_{t_{j-1}}^{\theta} - X_{s}^{\theta} \,
\dd W_s \Big\|_{L^{2} (\Omega ; \mathbb{R} )}
+ \Big\| \int_{t_{j-1}}^{t_{j}} \varphi(X_{t_{j-1}}) - \varphi(X_{s^{-}}) \,
\dd N_s \Big\|_{L^{2} (\Omega ; \mathbb{R} )}.
\end{aligned}
\end{equation}
As a result of Lemma \ref{lem:solution-estimation},
we know
\begin{equation}
\|X_{t_{j}}-X_{s}\|_{L^{2}(\Omega ; \mathbb{R})}
\leq C h^{\frac{1}{2}}.
\end{equation}
In addition, using the generalized It\^{o} formula \cite{Gardon2004approximations},
the Young inequality, the H\"older inequality, the It\^{o} isometry
and Lemmas \ref{lem:solution-moment-bound},
\ref{lem:solution-inverse-moment-bound}, we can obtain that
\begin{equation}\label{eq:solution-L2-esti}
\begin{aligned}
& \|X_{t_{j}}^{\gamma}-X_{s}^{\gamma}
\|^{2}_{L^{2}(\Omega ; \mathbb{R})}  \\
= & \Big\| \int_{s}^{t_{j}} \gamma \mu(X_r) X_r^{\gamma -1} \, \dd r
+ \int_{s}^{t_{j}} \gamma \phi(X_r) X_r^{\gamma -1} \, \dd W_r
+ \frac{1}{2}\gamma (\gamma-1)\int_{s}^{t_{j}}
X_{r}^{\gamma-2} \phi^{2}(X_r) \, \dd r \\
& \quad
+ \int_{s^{+}}^{t_{j}}\big[ (X_{r^{-}} + \varphi(X_{r^{-}}))^{\gamma}
- X_{r^{-}}^{\gamma} \big]  \dd N_{r} \Big\|^{2}_{L^{2}(\Omega ; \mathbb{R})}  \\
\leq &
4 \gamma^2 \Big\| \int_{s}^{t_{j}} \mu(X_r) X_r^{\gamma -1} \, \dd r
\Big\|^{2}_{L^{2}(\Omega ; \mathbb{R})}
+ 4 \gamma^2 \Big\| \int_{s}^{t_{j}} \phi(X_r) X_r^{\gamma -1} \, \dd W_r
\Big\|^{2}_{L^{2}(\Omega ; \mathbb{R})} \\
&   +
\gamma^2(\gamma-1)^2 \Big\| \int_{s}^{t_{j}} X_{r}^{\gamma-2}
\phi^{2}(X_r) \,\dd r \Big\|^{2}_{L^{2}(\Omega ; \mathbb{R})}
+
4  \Big\| \int_{s^{+}}^{t_{j}} \big[ \big(X_{r^{-}}
+ \varphi(X_{r^{-}})\big)^{\gamma} - X_{r^{-}}^{\gamma} \big]  \dd N_r
\Big\|^{2}_{L^{2}(\Omega ; \mathbb{R})}
\\
\leq &
4 \gamma^2 h  \int_{s}^{t_{j}} \E \big[|\mu(X_r) X_r^{\gamma -1}|^2\big] \dd r
+ 4 \gamma^2 \int_{s}^{t_{j}} \E \big[|\phi(X_r) X_r^{\gamma -1}|^2\big] \dd r \\
& +
\gamma^2(\gamma-1)^2 h \int_{s}^{t_{j}} \E\big[ | X_{r}^{\gamma-2}
\phi^{2}(X_r)|^2\big] \dd r
+
8 \lambda \int_{s^{+}}^{t_{j}} \E\big[ \big| \big(X_{r^{-}}
+ \varphi(X_{r^{-}})\big)^{\gamma} - X_{r^{-}}^{\gamma}\big|^2 \big]
\dd r \\
&  + 8\lambda^2 h \int_{s^{+}}^{t_{j}} \E \big[ \big| \big(X_{r^{-}}
+ \varphi(X_{r^{-}})\big)^{\gamma} - X_{r^{-}}^{\gamma}\big|^2 \big]
\dd r \\
\leq & C h .
\end{aligned}
\end{equation}
Similar to the proof of \eqref{eq:solution-L2-esti},
we can also derive that
\begin{equation}\label{eq:solution-inverse-L2-esti}
\begin{aligned}
& \|X_{t_{j}}^{-1}-X_{s}^{-1}
\|^{2}_{L^{2}(\Omega ; \mathbb{R})}  \\
\leq &
4  \Big\| \int_{s}^{t_{j}} \mu(X_r) X_r^{-2} \, \dd r
\Big\|^{2}_{L^{2}(\Omega ; \mathbb{R})}
+ 4  \Big\| \int_{s}^{t_{j}} \phi(X_r) X_r^{-2}\, \dd W_r
\Big\|^{2}_{L^{2}(\Omega ; \mathbb{R})} \\
&+
4 \Big\| \int_{s}^{t_{j}} X_{r}^{-3}
\phi^{2}(X_r) \, \dd r \Big\|^{2}_{L^{2}(\Omega ; \mathbb{R})}
+
4  \Big\| \int_{s^{+}}^{t_{j}} \big[ \big(X_{r^{-}}
+ \varphi(X_{r^{-}})\big)^{-1} - X_{r^{-}}^{-1} \big] \dd N_r
\Big\|^{2}_{L^{2}(\Omega ; \mathbb{R})}
\\
\leq &
4  h  \int_{s}^{t_{j}} \E [|\mu(X_r) X_r^{-2}|^2] \,\dd r
+ 4  \int_{s}^{t_{j}} \E [|\phi(X_r) X_r^{-2}|^2] \,\dd r
+
4 h \int_{s}^{t_{j}} \E[ | X_{r}^{-3}
\phi^{2}(X_r)|^2] \, \dd r \\
&+
8 \lambda \int_{s^{+}}^{t_{j}} \E\Big[
\Big| \tfrac{\varphi(X_{r^{-}})}
{ ( X_{r^{-}} + \varphi(X_{r^{-}})  )
	X_{r^{-}} }   \Big|^2 \Big]
\dd r  + 8\lambda^2 h \int_{s^{+}}^{t_{j}} \E\Big[
\Big| \tfrac{\varphi(X_{r^{-}})}
{ ( X_{r^{-}} + \varphi(X_{r^{-}})  )
	X_{r^{-}} }   \Big|^2 \Big] \dd r  \\
\leq &
4  h  \int_{s}^{t_{j}} \E [|\mu(X_r) X_r^{-2}|^2]\, \dd r
+ 4  \int_{s}^{t_{j}} \E [|\phi(X_r) X_r^{-2}|^2] \,\dd r
+
4 h \int_{s}^{t_{j}} \E[ | X_{r}^{-3}
\phi^{2}(X_r)|^2] \, \dd r \\
&+
8 \lambda(\lambda h + 1)\varepsilon_0^{-2}
\bigg(\int_{s^{+}}^{t_{j}} \E\Big[
\Big| \mathbbm{1}_{\{X_{r^{-}}<1\}}
\tfrac{\varphi(X_{r^{-}})}{X_{r^{-}}^{2}}
\Big|^2 \Big] \dd r
+
\int_{s^{+}}^{t_{j}} \E\Big[
\Big|  \mathbbm{1}_{\{X_{r^{-}}\geq 1\}}
\tfrac{\varphi(X_{r^{-}})}{X_{r^{-}}}
\Big|^2 \Big] \dd r \bigg)    \\
\leq &
4  h  \int_{s}^{t_{j}} \E [|\mu(X_r) X_r^{-2}|^2] \, \dd r
+ 4  \int_{s}^{t_{j}} \E [|\phi(X_r) X_r^{-2}|^2] \, \dd r
+
4 h \int_{s}^{t_{j}} \E[ | X_{r}^{-3}
\phi^{2}(X_r)|^2] \, \dd r \\
&+
C\int_{s^{+}}^{t_{j}} \Big(1 + \E\big[| {X_{r^{-}} |^{-2}}\big]
+
\E\big[| {X_{r^{-}} |^{-4}} \big]\Big)  \dd r   \\
\leq & C h .
\end{aligned}
\end{equation}
Gathering the above estimates together implies that
\begin{equation}\label{eq:drift-convergent-order}
\Big\|\int_{t_{j-1}}^{t_{j}} \mu\left(X_{t_{j}}\right)
-\mu\left(X_{s}\right) \, \dd s \Big\|_{L^{2}(\Omega ; \mathbb{R})}
\leq C h^{\frac{3}{2 }}.
\end{equation}
In a similar way together with the It\^{o} isometry,
one can show that
\begin{equation}\label{eq:Browian-convergent-order}
\Big\|\int_{t_{j-1}}^{t_{j}}X_{t_{j-1}}^{\theta}- X_{s}^{\theta} \,
\dd W_s \Big\|_{L^{2}(\Omega ; \mathbb{R})}^{2}
=\int_{t_{j-1}}^{t_{j}} \left\|X_{t_{j-1}}^{\theta}- X_{s}^{\theta}
\right\|_{L^{2}(\Omega ; \mathbb{R})}^{2} \mathrm{d} s
\leq  C h^2 .
\end{equation}
Employing the It\^o isometry, the H\"older inequality
and \eqref{Ass:jump-coefficient-Lipschitz} yields
\begin{equation}\label{eq:jump-esti}
\begin{aligned}
&  \Big\|\int_{t_{j-1}}^{t_{j}}
\varphi(X_{t_{j-1}})- \varphi(X_{s^{-}}) \,
\dd N_s \Big\|_{L^{2}(\Omega ; \mathbb{R})}^{2} \\
= &
\Big\|\int_{t_{j-1}}^{t_{j}}
\varphi(X_{t_{j-1}})- \varphi(X_{s^{-}}) \,
\dd ( \tilde{N}_s + \lambda s )
\Big\|_{L^{2}(\Omega ; \mathbb{R})}^{2}  \\
\leq &
2 \Big\|\int_{t_{j-1}}^{t_{j}}
\varphi(X_{t_{j-1}})- \varphi(X_{s^{-}}) \,
\dd \tilde{N}_s \Big\|_{L^{2}(\Omega ; \mathbb{R})}^{2}
+ 2 \lambda^2 \Big\|\int_{t_{j-1}}^{t_{j}}
\varphi(X_{t_{j-1}})- \varphi(X_{s^{-}}) \,
\dd s \Big\|_{L^{2}(\Omega ; \mathbb{R})}^{2} \\
\leq &
(2\lambda + 2 \lambda^2 h)\int_{t_{j-1}}^{t_{j}} \left\|
\varphi(X_{t_{j-1}})- \varphi(X_{s^{-}})
\right\|_{L^{2}(\Omega ; \mathbb{R})}^{2} \dd s \\
\leq &
C M^2 \int_{t_{j-1}}^{t_{j}}
\left\|
X_{t_{j-1}} - X_{s^{-}}\right\|_{L^{2}
	(\Omega ; \mathbb{R})}^{2} \dd s
\\
\leq & C h^2.
\end{aligned}
\end{equation}
In view of \eqref{eq:drift-convergent-order}, \eqref{eq:Browian-convergent-order}
and \eqref{eq:jump-esti}, one can arrive at
\begin{equation}\label{eq:remainder-term-convergence-order}
\left\|\mathcal{M}_{t_{j}}\right\|_{L^{2}(\Omega ; \mathbb{R})}
\leq C h .
\end{equation}
In order to bound 
$
\|\mathbb{E} ( \mathcal{M}_{t_{j}}
| \mathcal{F}_{t_{j-1}} ) \|_{L^{2}
	(\Omega ; \mathbb{R})}
$,
we first note that
\begin{equation}\label{eq:diffusion-condition-expectation}
\E \Big[ \int_{t_{j-1}}^{t_{j}}
b \big( X_{t_{j-1}}^{\theta}- X_{s}^{\theta} \big) \,
\dd W_s \Big| \mathcal{F}_{t_{j-1}}\Big] = 0 ,
\end{equation}
and use the compensated Poisson process and the Jensen type inequality to get
\begin{equation}\label{eq:poisson-jump-condition-expectation}
\begin{aligned}
&  \Big\| \E \Big(\int_{t_{j-1}}^{t_{j}}
\varphi(X_{t_{j-1}})- \varphi(X_{s^{-}})
\, \dd N_s
\Big | \mathcal{F}_{t_{j-1}}\Big) \Big\|_{L^{2}(\Omega ; \mathbb{R})}^{2} \\
\leq & 2\Big\| \E \Big(\int_{t_{j-1}}^{t_{j}}
\varphi(X_{t_{j-1}})- \varphi(X_{s^{-}})
\, \dd \tilde{N}_{s}
\Big | \mathcal{F}_{t_{j-1}}\Big) \Big\|_{L^{2}(\Omega ; \mathbb{R})}^{2} \\
& 
+ 
2\Big\| \lambda \E \Big(\int_{t_{j-1}}^{t_{j}}
\varphi(X_{t_{j-1}})- \varphi(X_{s^{-}})
\, \dd s
\Big | \mathcal{F}_{t_{j-1}}\Big) \Big\|_{L^{2}(\Omega ; \mathbb{R})}^{2} \\
\leq &
2 \lambda^2 h \int_{t_{j-1}}^{t_{j}} \left\|
\varphi(X_{t_{j-1}})- \varphi(X_{s^{-}})
\right\|_{L^{2}(\Omega ; \mathbb{R})}^{2} \dd s \\
\leq &
C M^2 h \int_{t_{j-1}}^{t_{j}}
\left\|
X_{t_{j-1}} - X_{s^{-}}\right\|_{L^{2}
	(\Omega ; \mathbb{R})}^{2}  \dd s
\\
\leq & C h^3.
\end{aligned}
\end{equation}
Furthermore, the Jensen inequality together with \eqref{eq:drift-convergent-order}
implies
\begin{equation}\label{eq:drift-condition-expectation-convergent-order}
\Big\|\E \Big(\int_{t_{j-1}}^{t_{j}} \mu (X_{t_{j}})
-\mu (X_{s}) \, \dd s | \mathcal{F}_{t_{j-1}} \Big)
\Big\|_{L^{2}(\Omega ; \mathbb{R})}
\leq \Big\|\int_{t_{j-1}}^{t_{j}} \mu (X_{t_{j}})
-\mu (X_{s}) \, \dd s
\Big\|_{L^{2}(\Omega ; \mathbb{R})}
\leq C h^{\frac{3}{2 }}.
\end{equation}
Armed with \eqref{eq:diffusion-condition-expectation},
\eqref{eq:poisson-jump-condition-expectation} and
\eqref{eq:drift-condition-expectation-convergent-order},
one can directly derive that
\begin{equation}\label{eq:remainder-term-condition-expectation-convergence-rate}
\left\|\mathbb{E}\left(\mathcal{M}_{t_{j}}
| \mathcal{F}_{t_{j-1}}\right)\right\|_{L^{2}
	(\Omega ; \mathbb{R})}
\leq C h^{\frac{3}{2 }}.
\end{equation}
Plugging \eqref{eq:remainder-term-convergence-order} and \eqref{eq:remainder-term-condition-expectation-convergence-rate} into \eqref{eq:error-bound} finishes the proof of
Theorem \ref{thm:convergence-rate}.
\qed
\subsection{The critical case $ \gamma+1=2 \theta $}
Next we turn to the error analysis of the considered scheme for the model in the critical case $ \gamma+1=2 \theta $. 
Similar to Lemma \ref{lem:solution-estimation}, we first get the next lemma.
\begin{lemma} \label{lem:solution-esti-critical-case}
	Let all conditions in Proposition \ref{prop:Ait-Sahalia-solution-well-posedness} hold with model parameters satisfying 
	$ \gamma + 1 = 2 \theta $ and $ \frac{a_2}{b^2} > 2 \gamma - \frac{3}{2} $.
	Let $ \{ X_t \}_{ t \geq 0 } $ be the process defined by \eqref{eq:solution-Ait-Sahalia-model}.
	Then it holds that
	\begin{equation}\label{eq:solution-L2-esti-critical-case}
	\left\|X_{t}-X_{s}\right\|_{L^{2}(\Omega ; \mathbb{R})}
	\leq C|t-s|^{\frac{1}{2}}.
	\end{equation}
\end{lemma}
\emph{Proof of Lemma \ref{lem:solution-esti-critical-case}.}
Similar to the proof of Lemma \ref{lem:solution-estimation},
one can show that
\begin{equation}\label{eq:solution-L2-bound-critical-case}
\begin{split}
\E [ | X_{t}-X_{s}|^2]
\leq &
4 |t-s|   \int_{s}^{t} \E[ | \mu(X_r) |^2] \, \dd r
+ 4   \int_{s}^{t}
\E [| \phi(X_{r})|^2 ] \,\dd r
+ 4   \int_{s}^{t} \E [|\varphi(X_{r^{-}})|^2]
\, \dd r \\
& + 4 \lambda^2 |t-s|
\int_{s}^{t} \E [| \varphi(X_{r^{-}}) |^2 ] \, \dd r .
\end{split}
\end{equation}
It is easy to check that
\begin{equation}\label{eq:condition-critical}
\tfrac{2a_2 + b^2}{b^2} = \tfrac{2 a_2 }{b^2} +1 >
4 \gamma - 2 > 2 \gamma > 2 \theta
\end{equation}
under the condition
$ \frac{a_2}{b^2} > 2 \gamma - \frac{3}{2} $.
With the aid of Lemmas \ref{lem:solution-moment-bound},
\ref{lem:solution-inverse-moment-bound}, we can further deduce that
\begin{equation}\label{eq:solution-error-order}
\begin{split}
\E [ | X_{t}-X_{s}|^2]
\leq &
C |t-s| \int_{s}^{t} \big(1+\E [ | X_r |^{-2} ]
+ \E [ | X_r |^{2 } ] +
\E [ | X_r |^{2\gamma } ] \big) \, \dd r \\
& + C  \int_{s}^{t}
\E [ | X_r |^{2 \theta} ] \, \dd r
+ C (1 + |t-s|) \int_{s}^{t}
\big( 1 + \E[ | X_{r^{-}}|^{2} ] \big) \, \dd r \\
\leq &
C |t-s|^2  + C | t -s |
+ C | t -s |
+ C | t -s |^2 \\
\leq & C | t -s |.
\end{split}
\end{equation}
The proof is finished.
\qed
\begin{theorem}\label{thm:convergence-rate-critical-case}
	Let all conditions in Proposition \ref{prop:Ait-Sahalia-solution-well-posedness} hold with model parameters satisfying 
	$ \gamma + 1 = 2 \theta $ and $ \frac{a_2}{b^2} > 2 \gamma - \frac{3}{2} $.
	Let $ \{ X_{t} \}_{0 \leq t \leq T } $
	and $ \{ Y_{n} \}_{0 \leq n \leq N } $ be solutions to \eqref{eq:Ait-Sahalia-model}
	and \eqref{eq:BEM-numerical-approximation}, respectively.
	Then for any
	$ h < \frac{1}{2a_1} $ we have
	\begin{equation}\label{eq:convergence-rate-critical-case}
	\sup _{N \in \mathbb{N}} \sup _{0 \leq n \leq N}
	\|Y_{n} -  X_{t_{n}} \|_{L^{2}
		(\Omega ; \mathbb{R})} \leq C h^{\frac{1}{2}}.
	\end{equation}
\end{theorem}
\emph{Proof of Theorem \ref{thm:convergence-rate-critical-case}.}
The proof here follows the same lines in the proof of 
Theorem \ref{thm:convergence-rate} for the case $ \gamma + 1 > 2 \theta $.
Clearly, equalities \eqref{eq:solution-numerical-difference}-\eqref{eq:En-expansion-expectation} still hold for the critical case $ \gamma + 1 = 2 \theta $ and
we need to estimate \eqref{eq:remainder-term-bound} first. Similar to \eqref{eq:proof-remainder-term-bound}, we have
\begin{equation}
\begin{split}
\E [| \mathcal{M}_{t_{n}} |^2]
\leq &
4 h  \int_{t_{n-1}}^{t_n}
\E \big[| \mu(X_{t_n}) - \mu(X_s)|^2\big]  \, \dd s
+ 4 \int_{t_{n-1}}^{t_n}
\E \big[ | \phi(X_{t_{n-1}}) - \phi(X_s)|^2\big] \, \dd s\\
&  +
4 \lambda \int_{t_{n-1}}^{t_n}
\E \big[| \varphi(X_{t_{n-1}}) - \varphi(X_{s^{-}})|^2 \big] \, \dd s
+ 4 \lambda^2 h \int_{t_{n-1}}^{t_n}
\E \big[| \varphi(X_{t_{n-1}}) - \varphi(X_{s^{-}})|^2\big] \, \dd s  \\
\leq &
24 h   \int_{t_{n-1}}^{t_n}
a_{-1}^{2} ( \E[| X_{t_n}^{-1} |^2] + \E[| X_{s}^{-1} |^2])
+ a_1^{2} ( \E[| X_{t_n} |^2] + \E[| X_{s} |^2]) \\
&
+ a_2^{2} ( \E[| X_{t_n} |^{2\gamma}] + \E[| X_{s} |^{2\gamma}])
\,
\dd s +
8 b^2\int_{t_{n-1}}^{t_n}
\E[| X_{t_{n-1}} |^{2\theta}] + \E[| X_{s} |^{2\theta}] \, \dd s \\
&
+ 8 \lambda M^2 (1 + \lambda h) \int_{t_{n-1}}^{t_n}
\E[| X_{t_{n-1}}|^{2}] +
\E[| X_{s^{-}}|^{2}]  \, \dd s.
\end{split}
\end{equation}
This together with \eqref{eq:condition-critical} and
Lemmas \ref{lem:solution-moment-bound}, \ref{lem:solution-inverse-moment-bound}
implies that
\begin{equation}
\E [| \mathcal{M}_{t_{n}} |^2]
 < \infty.
\end{equation}
Since
$ \frac{a_2}{b^2} > 2 \gamma - \frac{3}{2} $,
we can find some $ q > 2 $ such that
\begin{equation}
\begin{split}
\mu^{\prime}(x) + \tfrac{q-1}{2}|\phi'(x)|^{2}
& =
-a_{-1} x^{-2}+ a_{1}- a_{2} \gamma x^{\gamma-1}
+ \tfrac{q-1}{2} b^{2} \theta^{2}
x^{2 \theta-2}  \\
& <
a_{1} + \big( \tfrac{q-1}{2} b^{2} \theta^{2} -
a_{2} \gamma \big) x^{\gamma-1} \\
&   \leq a_{1}
, \quad x \in \mathbb{R}_{+},
\end{split}
\end{equation}
which implies
\begin{equation}\label{eq:couple-monotonicity-condition-critical-case}
\langle x-y, \mu(x)- \mu(y)\rangle + \tfrac{q-1}{2} |\phi(x)-\phi(y)|^{2}
\leq a_1 |x-y|^{2}, \quad \forall x, y \in {\R}_{+}.
\end{equation}
With this, one can also rely on the mathematical induction argument to acquire 
\eqref{eq:E_k-drift-bound} and \eqref{eq:diffusion-jump-bounds}
for  the critical case. Then repeating the same arguments used in  
\eqref{eq:convergen-error-order}-\eqref{eq:En-mean-square-expectation-esti}
one can arrive at \eqref{eq:error-bound} for the critical case.
So it remains to estimate
$ \E [ |\mathcal{M}_{t_{j}} |^{2} ] $ and
$ \E [ |\E (\mathcal{M}_{t_{j}} | \mathcal{F}_{t_{j-1}} ) |^{2} ] $, $ j = 1,2,\cdots,n $.
According to \eqref{eq:remainder-term-esti}, the estimate of $ \E [ |\mathcal{M}_{t_{j}} |^{2} ] $ relies on treating
$ \|X_{t_{j}}^{\gamma}-X_{s}^{\gamma}
\|^{2}_{L^{2}(\Omega ; \mathbb{R})} $,
$ \|X_{t_{j}}^{-1}-X_{s}^{-1}
\|^{2}_{L^{2}(\Omega ; \mathbb{R})}$
and
$\big\|\int_{t_{j-1}}^{t_{j}}X_{t_{j-1}}^{\theta}- X_{s}^{\theta}
\, \dd W_s \big\|_{L^{2}(\Omega ; \mathbb{R})}^{2}$
separately.
Owing to \eqref{Ass:jump-coefficient-lower-bound}, \eqref{eq:solution-L2-esti},
\eqref{eq:condition-critical}
and Lemmas \ref{lem:solution-moment-bound}, 
\ref{lem:solution-inverse-moment-bound},
we have
\begin{equation}\label{solution-gama-L2-esti-critical-case}
\begin{aligned}
& \|X_{t_{j}}^{\gamma}-X_{s}^{\gamma}
\|^{2}_{L^{2}(\Omega ; \mathbb{R})}  \\
\leq &
4 \gamma^2 h  \int_{s}^{t_{j}} \E \big[|\mu(X_r) X_r^{\gamma -1}|^2\big] \dd r
+ 4 \gamma^2 \int_{s}^{t_{j}} \E \big[|\phi(X_r) X_r^{\gamma -1}|^2\big] \dd r \\
& +
\gamma^2(\gamma-1)^2 h \int_{s}^{t_{j}} \E\big[ | X_{r}^{\gamma-2}
\phi^{2}(X_r)|^2\big] \dd r
+
8 \lambda(\lambda h +1) \int_{s^{+}}^{t_{j}} \E\big[ \big| \big(X_{r^{-}}
+ \varphi(X_{r^{-}})\big)^{\gamma} - X_{r^{-}}^{\gamma}\big|^2 \big]
\dd r \\
\leq &
C h  \int_{s}^{t_{j}} \big(\E \big[ | X_r |^{2\gamma -4} \big]
+ \E \big[ | X_r |^{2\gamma -2} \big] +
\E \big[ | X_r |^{2\gamma } \big] +
\E \big[ | X_r |^{4\gamma - 2} \big] \big) \dd r \\
& + 4 b^2 \gamma^2 \int_{s}^{t_{j}}
\E \big[ | X_r |^{3\gamma - 1} \big] \dd r
+ \gamma^2(\gamma-1)^2 b^2 h \int_{s}^{t_{j}}
\E\big[ | X_{r} |^{4\gamma-2}\big]  \dd r \\
& + C \int_{s^{+}}^{t_{j}}
\Big( \E\big[ | X_{r^{-}}|^{2\gamma} \big]
+ \E\big[ | \varphi(X_{r^{-}})|^{2\gamma} \big] \Big) \dd r \\
\leq &
C h  \int_{s}^{t_{j}} \big(\E [ | X_r |^{2\gamma -4} ]
+ \E \big[ | X_r |^{2\gamma -2} \big] +
\E \big[ | X_r |^{2\gamma } \big] +
\E \big[ | X_r |^{4\gamma - 2} \big] \big)\dd r \\
& + 4 b^2 \gamma^2 \int_{s}^{t_{j}}
\E \big[ | X_r |^{3\gamma - 1} \big] \dd r
+ \gamma^2(\gamma-1)^2 b^2 h \int_{s}^{t_{j}}
\E\big[ | X_{r} |^{4\gamma-2}\big] \dd r
+ C \int_{s^{+}}^{t_{j}}
\Big( 1 + \E\big[ | X_{r^{-}}|^{2\gamma} \big] \Big) \dd r \\
\leq & C h.
\end{aligned}
\end{equation}
In the same manner as above,
we derive from
\eqref{eq:solution-inverse-L2-esti} that
\begin{equation}\label{eq:solution-inverse-L2-esti-critical-case}
\begin{aligned}
& \|X_{t_{j}}^{-1}-X_{s}^{-1}
\|^{2}_{L^{2}(\Omega ; \mathbb{R})}  \\
\leq &
4  h  \int_{s}^{t_{j}} \E \big[|\mu(X_r) X_r^{-2}|^2\big] \dd r
+ 4  \int_{s}^{t_{j}} \E \big[|\phi(X_r) X_r^{-2}|^2\big] \dd r \\
& +
4 h \int_{s}^{t_{j}} \E\big[ | X_{r}^{-3}
\phi^{2}(X_r)|^2\big] \dd r
+
8 \lambda \int_{s^{+}}^{t_{j}} \E\Big[
\Big| \tfrac{\varphi(X_{r^{-}})}
{ ( X_{r^{-}} + \varphi(X_{r^{-}})  )
	X_{r^{-}} }   \Big|^2 \Big]
\dd r
\\
&  + 8\lambda^2 h \int_{s^{+}}^{t_{j}} \E\Big[
\Big| \tfrac{\varphi(X_{r^{-}})}
{ ( X_{r^{-}} + \varphi(X_{r^{-}})  )
	X_{r^{-}} }   \Big|^2 \Big] \dd r  \\
\leq &
C h \int_{s}^{t_{j}} \E \big[ | X_{r} |^{-6} \big]
+ \E \big[ | X_{r} |^{-4} \big] + \E \big[ | X_{r} |^{-2} \big]
+ \E \big[ | X_{r} |^{2\gamma - 4} \big] \dd r \\
& + C \int_{s}^{t_{j}}
\E \big[ | X_{r} |^{2\theta - 4} \big] \dd r
+ C \int_{s}^{t_{j}}
\E \big[ | X_{r} |^{4\theta - 6} \big] \dd r \\
& +
8 \lambda(\lambda h + 1)\varepsilon_0^{-2}
\bigg(\int_{s^{+}}^{t_{j}} \E\Big[
\Big| \mathbbm{1}_{\{X_{r^{-}}<1\}}
\tfrac{\varphi(X_{r^{-}})}{X_{r^{-}}^{2}}
\Big|^2 \Big] \dd r
+
\int_{s^{+}}^{t_{j}} \E\Big[
\Big|  \mathbbm{1}_{\{X_{r^{-}}\geq 1\}}
\tfrac{\varphi(X_{r^{-}})}{X_{r^{-}}}
\Big|^2 \Big] \dd r \bigg)    \\
\leq &
C h \int_{s}^{t_{j}} \E \big[ | X_{r} |^{-6} \big]
+ \E \big[ | X_{r} |^{-4} \big] + \E \big[ | X_{r} |^{-2} \big]
+ \E \big[ | X_{r} |^{2\gamma - 4} \big] \dd r \\
& + C \int_{s}^{t_{j}}
\E \big[ | X_{r} |^{2\theta - 4} \big] \dd r
+ C \int_{s}^{t_{j}}
\E \big[ | X_{r} |^{4\theta - 6} \big] \dd r
+ C \int_{s^{+}}^{t_{j}}
\Big(1 + \E \big[ | X_{r^{-}} |^{-2}
+ | X_{r^{-}} |^{-4} \big]\Big) \dd r  \\
\leq & C h .
\end{aligned}
\end{equation}
Similar to \eqref{solution-gama-L2-esti-critical-case}, we have
\begin{equation}
\|X_{t_{j}}^{\theta}-X_{s}^{\theta}
\|^{2}_{L^{2}(\Omega ; \mathbb{R})}
\leq C h.
\end{equation}
Therefore, for $ \gamma +1 = 2 \theta $ and
$ \frac{a_2}{b^2} > 2 \gamma - \frac{3}{2} $, it holds that
\begin{equation}
\|\mathcal{M}_{t_{j}}\|_{L^2(\Omega;\R)} \leq C h.
\end{equation}
Following similar arguments as used in \eqref{eq:remainder-term-condition-expectation-convergence-rate}
and considering $ \gamma +1 = 2 \theta $,
$ \frac{a_2}{b^2} > 2 \gamma - \frac{3}{2} $,
we can obtain
\begin{equation}\label{eq:remainder-term-condition-expectation-critical-case}
\left\|\mathbb{E}\left(\mathcal{M}_{t_{j}}
| \mathcal{F}_{t_{j-1}}\right)\right\|_{L^{2}
	(\Omega ; \mathbb{R})}
\leq C h^{\frac{3}{2 }}.
\end{equation}
Based on \eqref{eq:error-bound},
we can easily derive that
\begin{equation}
\sup _{N \in \mathbb{N}} \sup _{0 \leq n \leq N}
\|Y_{n} -  X_{t_{n}} \|_{L^{2}
	(\Omega ; \mathbb{R})} 
	\leq C h^{\frac{1}{2}},
\end{equation}
as required.
\qed
\section{Numerical results}
In the present section, we provide some numerical experiments to support
the previous findings. We consider a generalized Ait-sahalia-type
rate model with Poisson jumps of the form
\begin{equation}\label{eq:Ait-Sahalia-model-Numerical}
\dd X_t= (a_{-1} X_t^{-1} - a_{0}+ a_{1}
X_t - a_{2} X_t^{\gamma} ) \, \dd t
+b X_t^{\theta} \dd W_t + \varphi(X_{t^{-}}) \, \dd N_{t},
\quad X_0= 1 > 0,
\end{equation}
with $ a_{-1}, a_{0}, a_{1}, a_{2}, b >0 $
and $ \gamma ,\theta > 1 $.
In the following tests we take two sets of model parameters:
\begin{itemize}
	\item Case 1: $ a_{-1}=2,
	a_{0} = 1, a_{1} = 1.5, a_{2} = 5,b=1, \theta=2,\gamma=3.5,\lambda=1 $;
	\item Case 2: $ a_{-1}=2,
	a_{0} = 1, a_{1} = 1.5, a_{2} = 5,b=1, \theta=2,\gamma=3,\lambda=1 $.
\end{itemize}
The first set of parameters satisfies $ \gamma + 1 > 2\theta $
and the second one belongs to the critical case $ \gamma + 1 = 2\theta $.
As the first part of numerical results,  we list in
Table 1 the percentage of negative paths by using the explicit
Euler method (EM) and the BEM with three time stepsizes
$ h=\frac{1}{4}, \frac{1}{8},\frac{1}{16} $ over $ 10^5 $ paths.
It is clear that the EM method fails to preserve the positivity,
although the percentage of negative paths decreases as the stepsize becomes
smaller. However, as expected, the BEM remains positive for
all stepsize, which is consistent with the previous theoretical results.
\begin{table}[h]
	\begin{tabular}{*{6}{c}}
		\toprule
		Stepsize ($T=1$) & \quad $ \varphi(x) $ & \quad EM (Case 1)& \quad
		EM (Case 2) & \quad
		BEM (Case 1) & \quad
		BEM (Case 2) \\
		\midrule
		\multirow{3}*{$\dfrac{1}{4}$} & $ -0.2x $ & $ 89.91\% $
		& $ 86.75\% $ & $ 0 $ & $ 0 $ \\&
		$ x $ & $ 91.48\% $
		& $ 88.46\% $ & $ 0 $ & $ 0 $ \\ &
		$ \sin(x) $ & $ 91,09\% $
		& $ 88.03\% $ & $ 0 $ & $ 0 $  \\
		\midrule
		\multirow{3}*{$\dfrac{1}{8}$} & $ -0.2x $ & $ 47.81\% $
		& $ 40.52\% $ & $ 0 $ & $ 0 $ \\ &
		$ x $ & $ 71.51\% $
		& $ 64.31\% $ & $ 0 $ & $ 0 $ \\ &
		$ \sin(x) $ & $ 69.71\% $
		& $ 61.42\% $ & $0 $ & $ 0 $  \\
		\midrule
		\multirow{3}*{$\dfrac{1}{16}$} & $ -0.2x $ & $ 6.38\% $
		& $ 4.26\% $ & $ 0 $ & $ 0 $ \\ &
		$ x $ & $ 39.82\% $
		& $ 28.58\% $ & $ 0 $ & $ 0 $ \\ &
		$ \sin(x) $ & $ 33.35\% $
		& $ 22.91\% $ & $ 0 $ & $ 0 $  \\
		\bottomrule
	\end{tabular}
	\caption{The percentage of negative paths for EM and BEM with three
		stepsizes over $ 10^5 $ paths.}
\end{table}

\begin{figure}[tp]
	\includegraphics[width=0.5\linewidth, height=0.3\textheight]{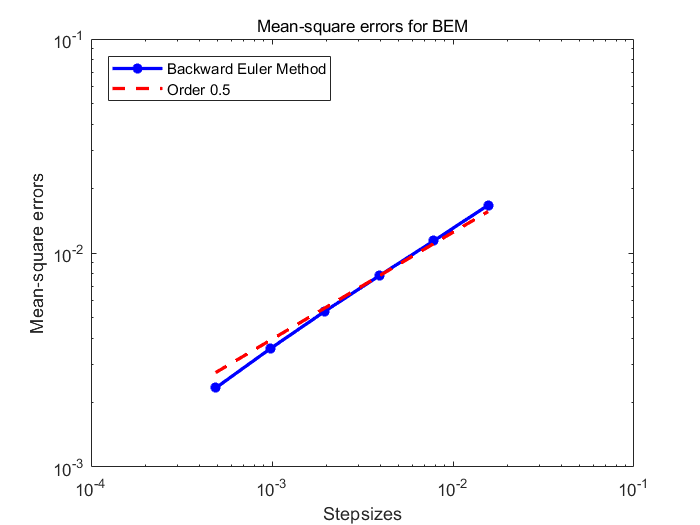}
	\includegraphics[width=0.5\linewidth, height=0.3\textheight]{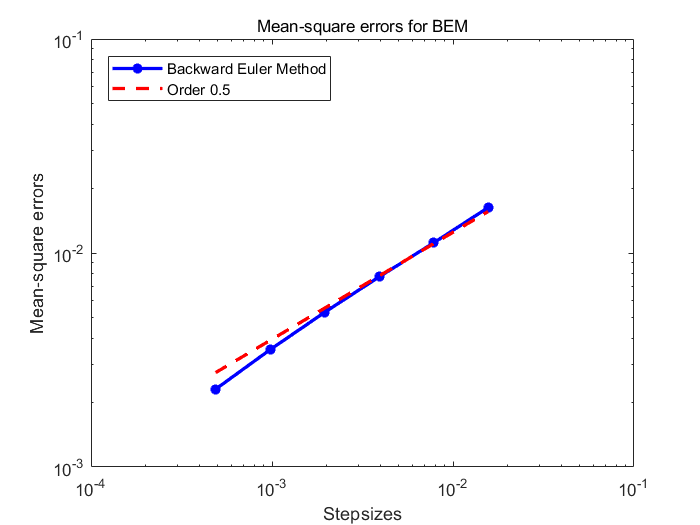}
	\caption[fig 1.]{Numerical results for \eqref{eq:Ait-Sahalia-model-Numerical}
		with $T=1:$ Mean-square convergence rate of BEM with
		$\varphi(x) = -0.2x $ for Case 1 (Left) and Case 2 (Right).}
	\label{fig:numerical-test_-0.2x}
\end{figure}

\begin{figure}[tp]
	\includegraphics[width=0.5\linewidth, height=0.3\textheight]{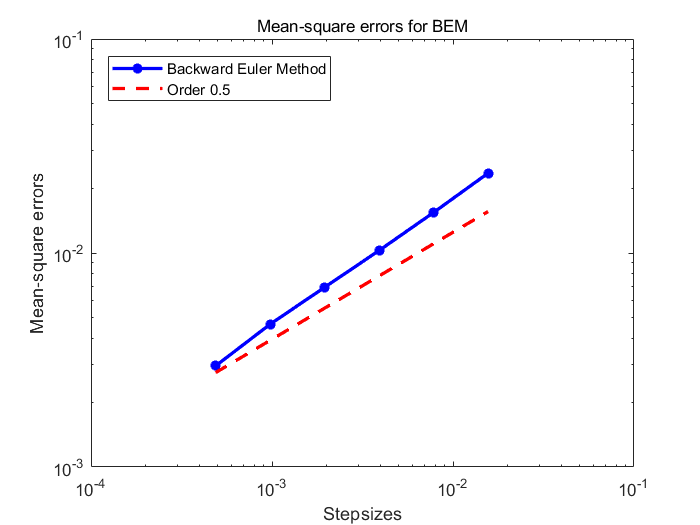}
	\includegraphics[width=0.5\linewidth, height=0.3\textheight]{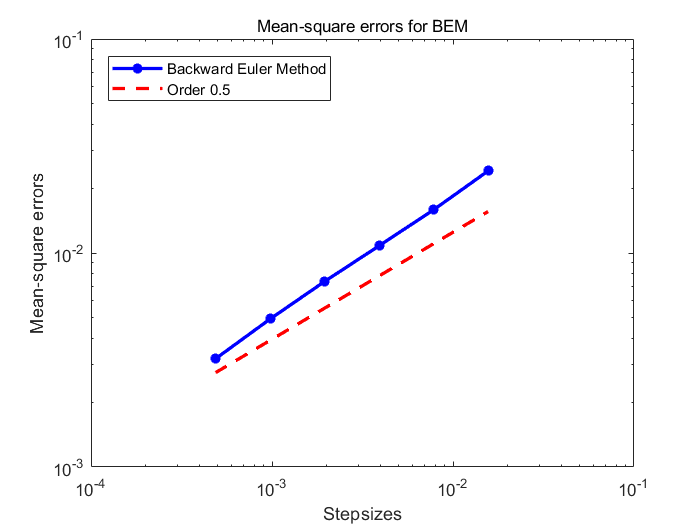}
	\caption[fig 1.]{Numerical results for \eqref{eq:Ait-Sahalia-model-Numerical}
		with $T=1:$  Mean-square convergence rate of BEM with
		$\varphi(x) = x $ for Case 1 (Left) and Case 2 (Right).}
	\label{fig:numerical-test_x}
\end{figure}

\begin{figure}[t]
	\includegraphics[width=0.5\linewidth, height=0.3\textheight]{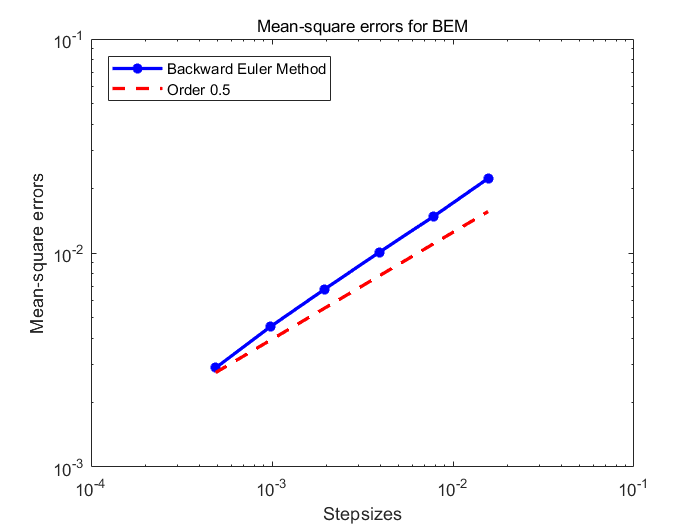}
	\includegraphics[width=0.5\linewidth, height=0.3\textheight]{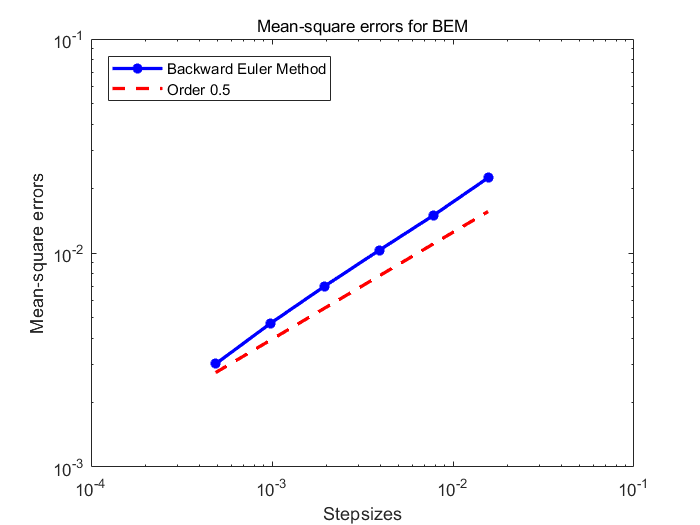}
	\caption[fig 1.]{Numerical results for \eqref{eq:Ait-Sahalia-model-Numerical}
		with $T=1:$  Mean-square convergence rate of BEM with
		$\varphi(x) = sin(x) $ for Case 1 (Left) and Case 2 (Right).}
	\label{fig:numerical-test_sinx}
\end{figure}

In the next step, we are to test the mean-square convergence
rate of BEM for six different stepsizes $ h = 2^{-j}, j=7,8,\cdots,11 $. As
usual, the expectation is approximated by using $ 10000 $ Brownian
and Poisson paths and the "true" solution of the model is identified
with the numerical one using a small stepsize $ h_{\text{exact}} = 2^{-13} $.
In order to clearly display the convergence rates, for various choices of the jump coefficient function $ \varphi $ we 
depict in Figures \ref{fig:numerical-test_-0.2x},
\ref{fig:numerical-test_x},
\ref{fig:numerical-test_sinx} mean-square approximation errors 
against six different stepsizes on a log-log scale.
From Figures \ref{fig:numerical-test_-0.2x},
\ref{fig:numerical-test_x},
\ref{fig:numerical-test_sinx}, one can observe that the mean-square error (solid lines)
and the reference (dashed lines) match well,
which indicates the mean-square convergence rate of order one-half.


	
\bibliographystyle{siam}
\bibliography{Ait_Sahalia_with_Poisson_Jump}
	
\end{document}